\documentclass[11pt]{article}
\usepackage{latexsym,euscript,amsmath,amssymb,amsbsy,amsfonts,amsthm,amsopn,amstext,amsxtra,amscd}
\usepackage{epsfig}
\usepackage{graphics}
\usepackage{url}

\usepackage[english]{babel}

\theoremstyle{plain}

\theoremstyle{definition}

\newcommand{\N}{{\mathbb N}}

\begin{document}
\centerline{\bf ON ONE TYPE OF COMPACTIFICATION}
\centerline{\bf OF POSITIVE INTEGERS}
\vskip0,5cm
\centerline{ Milan Pasteka}
\vskip0,5cm
Abstract.{\it The object of observation is a compact metric ring containing positive integers as dense subset. It is proved that this ring is isomorphic with 
a ring of reminder classes of polyadic integers.} 
$$ $$
Let $\N$ be the set of positive integers.
A mapping $|| \cdot  ||: \N \to <0,\infty)$ will be called norm if and only if the following conditions are satisfied for $a,b \in \N$ \newline
$|| a || = 0 \Leftrightarrow a=0$, \newline
$|| a+b  || \le || a  ||+|| b  ||$, \newline
$|| ab  || \le || b  ||.$\newline

There are various examples of norms on $\N$. One of these is polyadic norm defined in [N], [N1]. We start by a generalization of polyadic norm.
Denote by $a+(m)$ the arithmetic progression with difference $m$ which contains $a$. Instead of $0+(m)$ we write only $(m)$.

A subset $A \subset \N$ we call {\it closed to divisibility} or shortly CD - set if and only if \newline
$1 \in A,$ \newline
$m \in A, d|m \Rightarrow d \in A,$ \newline
$m_1, m_2 \in A \Rightarrow [m_1, m_2] \in A$, for $d,m,m_1,m_2 \in \N$.

Suppose that $A$ is infinite CD - set and $\{B_n\}$ is such sequence elements of $A$ that for every $d\in A$ there exists $n_0$ that $d|B_n$ for $n>n_0$.
It is easy to see that the mapping
$$
||a||_A = \sum_{n=1}^\infty \frac{h_n(a)}{2^n}
$$
for $a \in \N$ where $h_n(a) = 1-\mathcal{X}_{(B_n)}$, is a norm. This norm will be called {generalized polyadic} norm and the completion with respect the metric given by this norm will be called the  ring of generalized polyadic integers.

If $A=\N$ and $B_n = n!$ we get polyadic norm and the completion will be the ring of polyadic integers. In the case $A=\{ p^n; n=0,.,2,..\}$
and $B_n=p^n$ for given prime $p$ we obtain $p$ - adic norm and the completion will be the ring of $p$ - adic integers.

$$ $$

In the following text we shall assume that is given a compact metric space $(\Omega, \rho)$ containing $\N$ as dense subset. We suppose that the operations
 addition and multiplication on $\N$ are continuous and are extended to whole $\Omega$ to continuous operations. Thus $(\Omega, +, \cdot)$ is a topological
 commutative semiring.

Since $\Omega$ is compact we can suppose that there exists an increasing sequence of positive integers $\{x_n\}$ convergent to an element of $\Omega$.
Put $a_n = x_{2n} - x_n, n =1,2,...$. Then $a_n \ge n$ and
$$
a_n \to 0  \eqno{(1)}
$$
in the topology of $\Omega$.

For
$\beta \in {\Omega}$ and $b_n \to \beta$, $b_n \in \N$ we can consider the sequence of positive integer
$\{a_{k_n} - b_n \}$,  for suitable increasing sequence $\{k_n\}$, that
$a_{k_n} - b_n \to \beta'$ where $\beta + \beta'=0$. We see that $({\Omega}, +)$ is a compact group.

Clearly for every
$m \in \N$ there holds $cl(r+(m))= r+m{\Omega}$ where $m{\Omega}$ is the principal ideal in the ring
$({\Omega}, +, \cdot)$ generated by $m$. This yields
$$
{\Omega} = m{\Omega} \cup (1+m{\Omega}) \cup ...\cup (m-1+m{\Omega}). \eqno{(2)}
$$
Since the divisibility by $m$ in $\N$ is not necessary equivalent with divisibility by $m$ in
${\Omega}$ it is not provided that the last decomposition is disjoint.

{\bf Lemma 1.} {\bf Let $m \in \N$ is such positive integer that it is also the minimal
generator of the ideal $m {\Omega}$. Then every positive integer is divisible by $m$ in $\N$
if and only if it is divisible by $m$ in ${\Omega}$.}

{\bf Proof.} One implication is trivial. Suppose now that some positive integer $a$ is divisible
by $m$ in ${\Omega}$. Thus $a \in m{\Omega}$. Put $d=(a,m)$ - the greatest common divisor in $\N$. Then for $d = ax + my$ certain integers $x,y$.
This yields $d \in m{\Omega}$. We get $d{\Omega} = m{\Omega}$ and the minimality of
$m$ implies $m=d$. \qed

For every $ n \in \N$ we can define $g(n)$ as the minimal positive generator of $n{\Omega}$. Put
$\mathcal{A}=\{g(n); n \in \N\}$.

The set $r+m{\Omega}$ is closed and so from $(2)$ we see that also open, so called {\it clopen}.

It is easy to check that the set $\mathcal{A}$ is  CD - set.

Let $\{a_n\}$ be the sequence of positive integers given in (1). Clearly
$$
\bigcap_{n=1}^\infty a_n\Omega = \{0\},
$$
this yields
$$
\bigcap_{m \in \mathcal{A}} m\Omega = \{0\}. \eqno{(3)}
$$
And so we ge that the set $\mathcal{A}$ is infinite.  Since $m\Omega$ is open for $m \in \mathcal{A}$ equality  (3) implies that for each sequence $\{\alpha_n\}$ there holds
$$
\alpha_n \to 0 \Longleftrightarrow \forall m \in \mathcal{A} \exists n_0 ; n \ge n_0 \Longrightarrow m|\alpha_n.
$$
If we define the congruence  by the natural manner:  $\alpha \equiv \beta \pmod{\gamma}$ if and only if $\gamma$ divides $\alpha - \beta$, for
$\alpha, \beta, \gamma \in \Omega$, then there holds :
$$
\alpha_n \to \beta \Longleftrightarrow \forall m \in \mathcal{A} \exists n_0 ; n \ge n_0 \Longrightarrow \alpha_n \equiv \beta \pmod{m}.
$$
Thus the convergence can be metrised by the generalized polyadic norm : Let $\mathcal{A}=\{m_n, n=1,2,...\}$ and $M_n = [m_1,...,m_n], n=1,2,...$ then
$$
||\alpha||_{\mathcal{A}}=\sum_{n=1}^\infty \frac{1-\mathcal{X}_{M_n\Omega}(\alpha)}{2^{-n}}.
$$
We get

{\bf Theorem 1.} {\bf The metric $\rho$ is equivalent with the metric $\rho_{\mathcal{A}}$ where $\rho_{\mathcal{A}}(\alpha, \beta)=||\alpha -\beta||_{\mathcal{A}}$ for $\alpha, \beta \in \Omega$.}

And so every set $S \subset \Omega$ we have
$$
cl(S) = \bigcap_{n=1}^\infty (S+M_n\Omega) \eqno{(3)}.
$$
Denote by ${P}$ the Haar probability measure defined on $({\Omega}, +)$.  For $m \in \mathcal{A}$ the decomposition (2) is disjoint and so
${P}(r+m{\Omega}) = \frac{1}{m}$.
If we define the submeasure $\nu^\ast$ on the system of subsets of $\N$ as
$\nu^\ast(S)=P(cl(S))$ we get from (3) and upper semicontinuity of measure that
for each $S$
$$
\nu^\ast(S)=\lim_{n \to \infty} \frac{R(S:M_n)}{M_n}
$$
where $R(S:M_n)$ the number of elements of $S$ incongruent modulo $M_n$.

{\bf Theorem 2.} {\bf Let $\alpha, \beta \in \Omega$. There exist $\alpha_1, \beta_1 \in \Omega$ that the element
$\delta = \alpha_1\alpha + \beta_1\beta$ divides $\alpha$ and $\beta$.}

{\bf Proof.} Let $\{a_n\}, \{b_n\}$ be the sequences of positive integers that $a_n \to \alpha, b_n \to \beta$. Let
$d_n$ the greatest common divisor of $a_n, b_n, n=1,2,...$. Then $d_n = v_na_n + u_nb_n$ for some
$u_n, v_n$ - integers. The compactness of $\Omega$ provides that $u_{k_n} \to \alpha_1$
and $v_{k_n} \to \beta_1$ for suitable increasing sequence $\{k_n\}$. Put $\delta  = \alpha_1\alpha + \beta_1\beta$. We see that
$d_{k_n} \to \delta$. For $n=1,2,...$ we have $a_{k_n} = c_n d_{k_n}$. Since $\{c_n\}$ contains a convergent subsequence we get
that $\delta$ divides $\alpha$. Analogously can be derived that $\delta$ divides $\beta$. \qed

The element $\delta$ from Theorem 2 will be called the {\it greatest common divisor} of $\alpha, \beta$ and we shall write
$\delta \sim (\alpha, \beta)$.

{\bf Corollary 1.} {\bf If $p \in \mathcal{A}$ is a prime then for every $\alpha \in \Omega$ there holds $p$ divides $\alpha$ or
$(\alpha, p) \sim 1$.}

{\bf Proof.} If $p$ does not divide $\alpha$ then $\alpha \in \Omega \setminus p\Omega$. Consider a sequence of positive integers
$\{a_n\}$ which converges to $\alpha$. The set $\Omega \setminus p\Omega$ is open, thus we can suppose that $(a_n, p)=1$. This yields
$\ell_na_n + s_np=1$ for suitable integers $\ell_n, s_n$. Since $\Omega$ is compact space there exists an increasing sequence 
$\{k_n\}$ that $\ell_{k_n} \to \lambda, s_{k_n} \to \sigma$. And so $\lambda\alpha + \sigma p=1$. \qed

Analogously can be proved

{\bf Corollary 2.} {\bf An element $\alpha \in \Omega$ is invertible if and only if $(\alpha, p) \sim 1$ for every prime $p \in \mathcal{A}$.}

{\bf Lemma 2.} {\bf Each closed ideal in $\Omega$ is principal ideal.}

{\bf Proof.} Let $I \subset \Omega$ be closed ideal. Let $\alpha \in I$. Denote by $I_\alpha$ he set of all divisors of $\alpha$ belonging to $I$. The compactness of $\Omega$ yields that $I_\alpha$ is a closed set. From Lemma 2 we get that  $\alpha, \beta \in I$ there exists $\delta \in I$ that
$I_\delta \subset I_\alpha \cap I_\beta$. And so by induction can be proved that $I_\alpha, \alpha \in I$ is a centered system of closed sets. Thus its intersection in non empty, it contains some element $\gamma$. Then $I = \gamma\Omega$. \qed

$$ $$
In the sequel we denote $\Omega$ the ring of polyadic integers, thus completion of $\N$ with respect to norm $||\cdot||_{\N}$ and we suppose that
a infinite CD - set $A$ is given. The completion of $\N$ with respect to norm $||\cdot||_{A}$ we denote as $\Omega_A$.

Lemma 2 provides that $\cap_{a \in A} a\Omega = \alpha\Omega$ for suitable $\alpha \in \Omega$. We prove

{\bf Theorem 3.} {\bf The ring $\Omega_A$ is isomorphic with the factor ring $\Omega / \alpha\Omega$.}

{\bf Proof.} There holds \newline
If a sequence of positive integers is Cauchy's with respect to $||\cdot||_{\N}$ then it is Cauchy's with respect to $||\cdot||_A$ also. \newline
If $\{a_n\}, \{b_n\}$ are sequences of positive integer then $|||a_n - b_n|||_{\N} \to 0 \Longrightarrow |||a_n - b_n|||_A \to 0$.

Therefore we can define a mapping $F: \Omega \to \Omega_A$ in the following way: If $\beta \in \Omega, b_n \to \beta$ with respect $||\cdot||_{\N}$ then
$F(\beta)$ is the limit of $\{b_n\}$ with respect to $||\cdot||_A$. Clearly $F$ is a surjective morphism with kernel $\cap_{a \in A} a\Omega$ and the assertion follows. \qed

Katedra matematiky a informatiky, Pedagogicka Fakulta, Trnavska Univerzita \newline
Priemyselna 4, TRNAVA, SLOVAKIA \newline
e - mail: pasteka@mat.savba.sk

\end{document}